# A Closed Form Solution for the Normal Form and Zero Dynamics of a Class of Nonlinear Systems

Siamak Tafazoli[1]

*Abstract*—The normal form and zero dynamics are powerful tools useful in analysis and control of both linear and nonlinear systems. There are no simple closed form solutions to the general zero dynamics problem for nonlinear systems. A few algorithms exist for determining the zero dynamics, but none is straightforward and all are difficult to apply to large dimensional problems. A Closed form solution to the zero dynamics problem would motivate more usage of this powerful technique. The author presents here a simple *algebraic methodology* for the normal form and zero dynamics calculation of a class of nonlinear systems, mostly found in dynamical mechanical systems. The solution is in closed form so that application of the theorem presented is straight forward. As an illustration, the zero dynamics calculations for the complex dynamics of a flexible spacecraft is presented to demonstrate the simplicity and usefulness of the proposed closed form solution.

*Index Terms*— Control, Differential Geometry, Normal Form, Zero Dynamics, Nonlinear Systems, Feedback Linearization, Flexible Spacecraft.

## I. Introduction

In this section, a brief review on normal forms and zero dynamics including an overview of the existing algorithms devised to obtain these characteristics are presented. Zero dynamics has been derived associated with the problem of zeroing the output while linearizing a nonlinear system locally. The notion of zero dynamics was introduced by [1]. By definition, the nonlinear system is locally transformed to the normal form and zero dynamics is defined as a differential equation in terms of parts of the transformed states. In other words, when the normal form is obtained, the zero dynamics can be determined. Thus, [1] provides the first algorithm to calculate zero dynamics by first calculating the normal form.

The differential geometric approach for calculating zero dynamics was suggested by [2]. Reference [3] has shown how to use [2] to derive the normal forms and hence obtain the zero dynamics. The interested reader is referred to [4], which provides a detailed account of the normal form and the zero dynamics concepts.

Reference [5] discusses the use of symbolic computation for finding the zero dynamics. For a "low" dimensional class of systems, the computation can be performed without human aid or intervention, making the zero dynamics procedure a feasible and valuable addition to the toolbox of the modeler, analyst, or control system designer. For system models that are more than moderately complex and of "high" dimension, symbolic computation cannot be fully enjoyed due to the complexity of parts of the algorithms that are (double) exponential in some measure of the problem size, or due to expression swell that cannot be easily eliminated. This implies that symbolic computation will not replace other tools, like those based on numerics, but will complement them.

[1] The author is with the Department of Electrical and Computer Engineering, Concordia University, Montreal, Canada, H3G 1M8 (Mak.Tafazoli@gmail.com)

Reference [6] extends the zero dynamics algorithm for nonlinear affine control systems to nonlinear non-affine control systems by using non-smooth analysis and viability theory. Reference [7] proposes to use differential algebra, in particular the Ritt algorithm to calculate zero dynamics of polynomial systems and shows that for affine polynomial Single-Input Single-Output (SISO) systems the Ritt algorithm gives a result which is equivalent to the ordinary zero dynamics. Furthermore, it is shown how the algorithm can be used to calculate zero dynamics for more general affine systems. It also discusses zero dynamics of generalized state-space descriptions. Reference [8] considers the problem for affine polynomial Multi-Input Multi-Output (MIMO) systems when the system has a vector relative degree. It also discusses how the algorithm can be used when the system is not polynomial or it has no vector relative degree. Reference [9] gives a definition of a normal form for systems on generalized state space form. Reference [10] shows how a normal form, corresponding to that of affine state space system, can be calculated for generalized state space descriptions. The calculations are performed by using Grobner bases.

A zero dynamics algorithm for general nonlinear control systems is proposed in [11], in which only very mild regularity assumptions are needed, namely, the constancy of the dimensions of certain distributions (and/or of the ranks of certain mapping) around a given point. From the zero dynamics algorithm the normal form for a general nonlinear control systems is presented, which is an extension of the affine nonlinear systems. Finally, one of the applications of the normal form to the problem of local feedback asymptotic stabilization is studied. This result does not require asymptotic stability in the first approximation for the zero dynamics, so that it may be useful in order to solve critical problems of local feedback asymptotic stabilization.

In [12] and [13], a method is proposed to derive the zero dynamics of physical systems from bond graph models. This method incorporates the definition of zero dynamics in the differential geometric approach and the causality manipulation in the bond graph representation. By doing so, the state equations of the zero dynamics can be obtained and the system elements which are responsible for the zero dynamics can be identified. Reference [12] deals with SISO systems and [13] deals with MIMO systems.

In [14], the concept of zero dynamics is introduced to nonlinear singular systems, and a zero dynamics algorithm is proposed. Normal form of nonlinear singular systems is obtained by using this zero dynamics algorithm. The decoupling problem is also studied via this normal form.

Reference [15] defines a generalized normal form and a computational algorithm is provided to obtain it. Reference [16] defines a generalized zero dynamics, introduces a generalized zero dynamics algorithm, and gives a sufficient condition which ensures stabilizability of linear systems with asymptotically stable generalized zero dynamics by static output feedback.

The above brief literature review demonstrates that, to date, different approaches to the problem of solving the zero dynamics problem are algorithmic/ computational using different mathematical approaches (i.e. differential geometry, differential algebra, symbolic computation, non-smooth analysis and viability theory, and bond graph models).

In what follows, to the best of the author's knowledge, he presents the first linear algebraic approach to the problem of finding the normal form and zero dynamics of a class of nonlinear systems which results in a *closed form solution*.

## II. MAIN RESULTS

***Theorem:*** Given a nonlinear system of the form

$$\dot{\vec{x}} = \vec{f}(\vec{x}) + \underline{G}(\vec{x})\,\vec{\tau} = \vec{f} + \sum_{j=1}^{p} \vec{g}_j \tau_j \qquad (1)$$
$$\vec{y} = \vec{h}(\vec{x})$$

where $\vec{x} \in \mathbb{R}^n$, $\vec{\tau} \in \mathbb{R}^p$ and $\vec{y} \in \mathbb{R}^s$ are the system state, input and output, respectively, and $\vec{f} \in \mathbb{R}^n$ and $\vec{h} \in \mathbb{R}^s$ are smooth vector fields belonging to class $C^\infty$. Furthermore, assume that the following conditions are satisfied:

(i) Distribution $\Gamma$, spanned by smooth vector fields $\vec{g}_j$, $j = 1, \ldots, p$ is involutive.

(ii) $p \leq s \leq r < n$, where $r$ is the well-defined relative degree of the system.

(iii) $\vec{f}$ can be expressed as $\underline{M}^{-1}(\vec{x}_\beta)\,\vec{l}(\vec{x})$, where $\underline{M} \in \mathbb{R}^{n \times n}$ is the invertible system matrix which is a function of only the state variables $\vec{x}_\beta = \begin{bmatrix} x_{p+1} & \cdots & x_n \end{bmatrix}^T$ and $\vec{l} \in \mathbb{R}^n$ contains all linear and non-linear terms.

Under the above conditions, the new *normal coordinates* $\vec{\eta} \in \mathbb{R}^{n-p}$, corresponding to $\vec{x}_\beta$ and given by the diffeomorphisms $\Phi_k$, $k = 1, \ldots, n-p$, is given by

$$\vec{\eta} = x_\alpha \underline{M}_{12} \underline{M}_{22}^{-1} + x_\beta \qquad (2)$$

where $\vec{x}_\alpha = \begin{bmatrix} x_1 & \cdots & x_p \end{bmatrix}^T$, and $\underline{M}_{12}$ and $\underline{M}_{22}$ are block matrices of $\underline{M}$ given by the following partitioning

$$\underline{M} = \begin{bmatrix} \underline{M}_{11} & \underline{M}_{12} \\ \underline{M}_{21} & \underline{M}_{22} \end{bmatrix} \qquad (3)$$

The rate of change of $\vec{\eta}$, which constitutes the key part of the *normal form* of the original nonlinear system (1), is given by

$$\dot{\vec{\eta}} = \underline{M}_{22}^{-1} \vec{l}_\beta + \left( \underline{M}_{22}^{-1} \underline{\dot{M}}_{12}^T + \underline{\dot{M}}_{22}^{-1} \underline{M}_{12}^T \right) \vec{\zeta} \qquad (4)$$

where $\vec{l}_\beta$ refers to the subset of vector $\vec{l}$ composed of $(p+1)^{th}$ to $n^{th}$ elements and $\vec{\zeta}$ is the new *normal coordinates* [4] corresponding to $\vec{x}_\alpha$ obtained using trivial mappings.

**Remark:** The matrix $\underline{G} \in \mathbb{R}^{n \times p}$ is composed of columns that are each smooth vector fields $\vec{g}_j$, $j = 1, \ldots, p$ defined on $\mathbb{R}^n$. It should be noted that the first two assumptions are quite standard and usual assumptions commonly found in the literature.

*Proof:* Let's start with (1) where the state variable $\vec{x} \in \mathbb{R}^n$ is partitioned into two parts as

$$\begin{aligned} \vec{x} &= \begin{bmatrix} x_1 & \cdots & x_p & | & x_{p+1} & \cdots & x_r & \cdots & x_n \end{bmatrix}^T \\ &= \begin{bmatrix} \vec{x}_\alpha^T & | & \vec{x}_\beta^T \end{bmatrix}^T \end{aligned} \qquad (5)$$

Let's assume that the $\vec{f}$ and $\underline{G}$ terms in (1) can be given as

$$\begin{aligned} \vec{f} &= \underline{M}^{-1}(\vec{x}_\beta) \vec{l}(\vec{x}) \\ \underline{G} &= \underline{M}^{-1}_{\{n,p\}}(\vec{x}_\beta) \end{aligned} \qquad (6)$$

where $\underline{M}$ is the system matrix (e.g. mass or inertial matrix in case of a mechanical system) and $\vec{l}$ is the combination of all linear and non-linear terms. The matrix $\underline{M}^{-1}_{\{n,p\}}$ is obtained by taking the first $p$ columns of matrix $\underline{M}^{-1}$. It has been assumed that $\underline{M}$ is a function of a subset of the state variables (i.e. $\vec{x}_\beta$) which is a reasonable assumption for many classes of physical systems (e.g. mechanical systems) such as a flexible spacecraft as will be described in the next section.

For most physical systems, we also have a symmetric matrix $\underline{M}$ such that $\underline{M}_{21} = \underline{M}_{12}^T$. The block matrices are specified as $\underline{M}_{11} \in \mathbb{R}^{p \times p}$, $\underline{M}_{12} \in \mathbb{R}^{p \times (n-p)}$, $\underline{M}_{21} \in \mathbb{R}^{(n-p) \times p}$ and $\underline{M}_{22} \in \mathbb{R}^{(n-p) \times (n-p)}$. Using the block matrix inversion formula [17], we can obtain

$$\underline{M}^{-1} = \begin{bmatrix} \underline{F}_{11}^{-1} & -\underline{F}_{11}^{-1} \underline{M}_{12} \underline{M}_{22}^{-1} \\ -\underline{M}_{22}^{-1} \underline{M}_{12}^T \underline{F}_{11}^{-1} & \underline{M}_{22}^{-1} + \underline{M}_{22}^{-1} \underline{M}_{12}^T \underline{F}_{11}^{-1} \underline{M}_{12} \underline{M}_{22}^{-1} \end{bmatrix} \qquad (7)$$

where $\underline{F}_{11} = \underline{M}_{11} - \underline{M}_{12} \underline{M}_{22}^{-1} \underline{M}_{12}^T$. For the sake of our subsequent development, the matrices $\underline{M}_{11}$, $\underline{M}_{22}$ and $\underline{F}_{11}$ must be nonsingular so that $\underline{M}^{-1}$ is invertible. The definition and specification of zero dynamics require the normal form of system (1) which can be obtained easily, once we have found the appropriate mappings (i.e. diffeomorphisms). In general, there are $n - p$ diffeomorphisms, $\Phi_k$, that can be solved for, of which only $n - r$ are non-trivial and of interest. All $\Phi_k$'s must respect the following conditions:

$$\begin{aligned} L_{\vec{g}_j} \Phi_k &= \frac{\partial \Phi_k}{\partial \vec{x}} \vec{g}_j = 0, \; k = 1, 2, \ldots, n-r; \; j = 1, \ldots, p \\ \Phi_k(0) &= 0 \\ \nabla \Phi_k & \text{ are linearly independent} \end{aligned} \qquad (8)$$

It's the first condition in (8) that can be used advantageously to formulate the problem of finding $\Phi_k$ in a new light. Let's express $L_{\vec{g}_j}\Phi_k = 0$ in its matrix from

$$L_{\vec{g}_j}\Phi_k = \nabla\Phi_k \vec{g}_j = 0, \ j = 1,\ldots,p \quad , \ k = 1,\ldots,n-p$$

$$\begin{bmatrix} \dfrac{\partial \Phi_k}{\partial x_1} & \cdots & \dfrac{\partial \Phi_k}{\partial x_n} \end{bmatrix} \begin{bmatrix} | & & | \\ \vec{g}_1 & \cdots & \vec{g}_p \\ | & & | \end{bmatrix} = \vec{0}^T \tag{9}$$

and now take the transpose of (9) to obtain

$$\begin{bmatrix} - & \vec{g}_1^T & - \\ & \vdots & \\ - & \vec{g}_p^T & - \end{bmatrix} \begin{bmatrix} \dfrac{\partial \Phi_k}{\partial x_1} \\ \vdots \\ \dfrac{\partial \Phi_k}{\partial x_n} \end{bmatrix} = \vec{0} \quad or \quad \underline{G}^T (\nabla\Phi_k)^T = \vec{0} \tag{10}$$

Hence, $(\nabla\Phi_k)^T = \mathbb{N}(\underline{G}^T)$ where $\mathbb{N}(\underline{G}^T)$ is the null space of $\underline{G}^T$. In essence, we can now partly solve the problem by finding the null space of $\underline{G}^T$ using linear algebraic techniques. The matrix $\underline{G}^T$ is the first $p$ rows of $\underline{M}^{-1}$ obtained from (7) as

$$\underline{G}^T = \begin{bmatrix} \underline{F}_{11}^{-1} & -\underline{F}_{11}^{-1}\underline{M}_{12}\underline{M}_{22}^{-1} \end{bmatrix} \tag{11}$$

Now assume that $\mathbb{N}(\underline{G}^T) = \underline{X}$ so that

$$\underline{G}^T \underline{X} = \underline{0} \tag{12}$$

and let $\underline{X} = \begin{bmatrix} \underline{N} \\ \underline{I} \end{bmatrix}$, where $\underline{I}$ is the identity matrix. Using (11) and (12), we solve for $\underline{X} \in \mathbb{R}^{n\times(n-p)}$

$$\underline{G}^T \underline{X} = \begin{bmatrix} \underline{F}_{11}^{-1} & -\underline{F}_{11}^{-1}\underline{M}_{12}\underline{M}_{22}^{-1} \end{bmatrix} \begin{bmatrix} \underline{N} \\ \underline{I} \end{bmatrix} = \underline{0} \tag{13}$$

Expanding (13) gives

$$\underline{F}_{11}^{-1}\underline{N} - \underline{F}_{11}^{-1}\underline{M}_{12}\underline{M}_{22}^{-1}\underline{I} = \underline{0} \tag{14}$$

By inspecting (14), we find that

$$\underline{N} = \underline{M}_{12}\underline{M}_{22}^{-1} \tag{15}$$

and hence

$$\mathbb{N}(\underline{G}^T) = \underline{X} = \begin{bmatrix} M_{12} M_{22}^{-1} \\ I \end{bmatrix} \quad (16)$$

The rank of $\underline{X}$ provides the number of independent solutions for $\Phi_k$, and as $\underline{X}$ includes the $(n-p) \times (n-p)$ identity matrix then there are no zero columns in $\underline{X}$ and $Rank(\underline{X}) = n - p$. From (16) we can solve for $\Phi_k, k = 1, \ldots, n-p$, which will be used to transform system (1) into its *normal form*. Each column in matrix (16) represents the gradient of $\Phi_k$ (i.e. $(\nabla \Phi_k)^T$) and expanding (16) gives

$$\begin{bmatrix} \dfrac{\partial \Phi_k}{\partial x_1} \\ \vdots \\ \dfrac{\partial \Phi_k}{\partial x_n} \end{bmatrix} = \begin{bmatrix} \overbrace{\begin{matrix} N_{11} & \cdots & N_{1(n-p)} \\ \vdots & \ddots & \vdots \\ N_{p1} & \cdots & N_{p(n-p)} \end{matrix}}^{p \times (n-p)} \\ -- \quad -- \quad --- \\ \underbrace{\begin{matrix} 1 & \cdots & 0 \\ \vdots & \ddots & \vdots \\ 0 & \cdots & 1 \end{matrix}}_{(n-p) \times (n-p)} \end{bmatrix}_{n \times (n-p)} \quad (17)$$

Recalling that $N_{jk}$ terms in (17) are not a function of $\vec{x}_\alpha = \begin{bmatrix} x_1 & \cdots & x_p \end{bmatrix}$, we can obtain the closed-form expression for $\Phi_k$ as

$$\Phi_k = x_1 N_{1k} + x_2 N_{2k} + \cdots + x_p N_{pk} + x_{p+k}.1 + c_k \quad \text{for} \quad 1 \leq k \leq n-p, \quad c_k \in \mathbb{R} \quad (18)$$

Let's define a vector containing all the $n - p$ diffeomorphisms $\Phi_k$ as $\vec{\Phi}^T = \begin{bmatrix} \Phi_1 & \Phi_2 & \cdots & \Phi_{n-p} \end{bmatrix}$, we can then express the analytical solution for all $\Phi_k$, which are in essence the new *normal variables* $\vec{\eta}$, as

$$\vec{\eta}^T = \vec{\Phi}^T = \begin{bmatrix} x_1 & \cdots & x_p \end{bmatrix} \underline{N} + \begin{bmatrix} x_{p+1} & \cdots & x_n \end{bmatrix} + \vec{c}^T \quad (19)$$

or simply as

$$\vec{\eta}^T = \begin{bmatrix} x_\alpha & | & x_\beta \end{bmatrix} \begin{bmatrix} M_{12} M_{22}^{-1} \\ I \end{bmatrix} = \vec{x}^T \underline{X} + \vec{c}^T \quad (20)$$

where $\vec{c} \in \mathbb{R}^{(n-p) \times 1}$ is a vector of arbitrary constants, which could be zero.

Having obtained $\vec{\eta} = \vec{\Phi}(\vec{x})$, we can solve for $\vec{x} = \vec{\Phi}^{-1}(\vec{\eta})$, knowing that the inverse map, $\vec{\Phi}^{-1}$, exists by definition. Once this is done, both $\underline{M}$ and $\vec{l}$ can be expressed as a function of the new *normal coordinates* $\vec{\eta}$ [4], that is, $\underline{M}(\vec{\eta})$ and $\vec{l}(\vec{\zeta}, \vec{\eta})$ where $\vec{\zeta} = \Psi(\vec{x})$ are the transformed $\vec{x}_\alpha$ states.

The corresponding normal form of (1) can then be obtained by differentiating (20), that is

$$\dot{\vec{\eta}}^T = \dot{\vec{x}}^T \underline{X} + \vec{x}^T \underline{\dot{X}} \tag{21}$$

Using equations (1), (6), (12) and (16) with (21), we can expand (21) as follows

$$\begin{aligned}\dot{\vec{\eta}}^T &= \left(\underline{M}^{-1} \vec{l} + \underline{G} \vec{\tau}\right)^T \underline{X} + \vec{x}^T \begin{bmatrix} \underline{\dot{N}} \\ 0 \end{bmatrix} \\ &= \vec{\tau}^T \underbrace{\underline{G}^T \underline{X}}_{=0} + \vec{l}^T \underline{M}^{-1} \underline{X} + \vec{x}_\alpha^T \underline{\dot{N}}\end{aligned} \tag{22}$$

Substituting (7) and (16) into (22) and simplifying, we obtain

$$\begin{aligned}\dot{\vec{\eta}}^T &= \vec{l}^T \begin{bmatrix} \underline{F}_{11}^{-1} & -\underline{F}_{11}^{-1} \underline{M}_{12} \underline{M}_{22}^{-1} \\ -\underline{M}_{22}^{-1} \underline{M}_{12}^T \underline{F}_{11}^{-1} & \underline{M}_{22}^{-1} + \underline{M}_{22}^{-1} \underline{M}_{12}^T \underline{F}_{11}^{-1} \underline{M}_{12} \underline{M}_{22}^{-1} \end{bmatrix} \\ &\quad \begin{bmatrix} \underline{M}_{12} \underline{M}_{22}^{-1} \\ \underline{I} \end{bmatrix} + \vec{x}_\alpha^T \underline{\dot{N}} \\ &= \vec{l}^T \begin{bmatrix} 0 \\ \underline{M}_{22}^{-1} \end{bmatrix} + \vec{x}_\alpha^T \left(\underline{\dot{M}}_{12} \underline{M}_{22}^{-1} + \underline{M}_{12} \underline{\dot{M}}_{22}^{-1}\right)\end{aligned} \tag{23}$$

The final compact form of (23) as a function of the new states $\vec{\zeta}$ and $\vec{\eta}$ is simply

$$\dot{\vec{\eta}} = \underline{M}_{22}^{-1}(\vec{\eta}) \vec{l}_\beta(\vec{\zeta},\vec{\eta}) + \left(\underline{M}_{22}^{-1}(\vec{\eta}) \underline{\dot{M}}_{12}^T(\vec{\eta}) + \underline{\dot{M}}_{22}^{-1}(\vec{\eta}) \underline{M}_{12}^T(\vec{\eta})\right) \vec{\zeta} \tag{24}$$

where $\vec{l}_\beta$ refers to the subset of vector $\vec{l}$ composed of $(p+1)^{th}$ to the $n^{th}$ elements.

***Corollary:*** The zero dynamics of system (1) is simply the normal equation (24) with all outputs kept at zero. Assuming that the output $\vec{y}$ in (1) is a linear function of the first $s$ elements of the state vector, $\vec{y} = \underline{C} \vec{x}_{\{1..s\}}$, where $\underline{C} \in \mathbb{R}^{s \times s}$ is any non-singular, constant real matrix, then the zero dynamics will have the simple form

$$\dot{\vec{\eta}} = \underline{M}_{22}^{-1} \vec{l}_\beta \tag{25}$$

In order to define the above, we need the non-singularity assumption since otherwise, $\underline{C}$ could have a null space and hence $\vec{x}_{\{1..s\}}$ need not be zero when $\vec{y} = 0$. Let's also assume that the new normal coordinates $\vec{\zeta}$ remains the same as $\vec{x}_{\{1..s\}}$, after the trivial mappings $\Psi(\vec{x})$. Note that $\vec{x}_{\{1..s\}}$ refers to the sub-vector of $\vec{x}$ composed of its first $s$ elements.

# III. AN ILLUSTRATIVE EXAMPLE: INTERNAL AND ZERO DYNAMICS OF A FLEXIBLE SPACECRAFT'S ATTITUDE CONTROL

The results presented in this paper are applicable to problems from different fields (e.g. robotics, general dynamics, etc.) as long as the stated assumptions are met. A particular field of interest is control of flexible robots and spacecraft having a nonlinear system dynamics which can be formulated in the representation (1). As a specific example of the use of the theorem presented in the previous section, we refer the interested reader to a MIMO system presented in [18] where the zero dynamics is key to the closed-loop stability analysis. The system is comprised of an input-output feedback linearizing controller applied to a spacecraft with flexible appendages and since this system is under-actuated, an internal/zero dynamics is present which needs to be analyzed.

The derivation of the system dynamics is not detailed here (refer to [19] for a detailed derivation), and the author presents only the resulting attitude and vibrational equations for a flexible spacecraft given by

$$\begin{bmatrix} \dot{\vec{\omega}} \\ \ddot{\vec{\chi}} \end{bmatrix} = \begin{bmatrix} \underline{F}_{11}^{-1} & -\underline{I}_t^{-1} \dot{\underline{\kappa}}_{\ddot{\chi}}^T \underline{F}_{22}^{-1} \\ -\underline{F}_{22}^{-1} \dot{\underline{\kappa}}_{\ddot{\chi}} \underline{I}_t^{-1} & \underline{F}_{22}^{-1} \end{bmatrix} \begin{bmatrix} \vec{l}_{\vec{\omega}} \\ \vec{l}_{\vec{\chi}} \end{bmatrix} + \begin{bmatrix} \underline{F}_{11}^{-1} \\ -\underline{F}_{22}^{-1} \dot{\underline{\kappa}}_{\ddot{\chi}} \underline{I}_t^{-1} \end{bmatrix} \vec{\tau} \tag{26}$$

where

$$\begin{aligned} \underline{F}_{11} &= \underline{I}_t - \dot{\underline{\kappa}}_{\ddot{\chi}}^T \underline{M}_{\phi\psi}^{-1} \dot{\underline{\kappa}}_{\ddot{\chi}} \\ \underline{F}_{22} &= \underline{M}_{\phi\psi} - \dot{\underline{\kappa}}_{\ddot{\chi}} \underline{I}_t^{-1} \dot{\underline{\kappa}}_{\ddot{\chi}}^T \\ \vec{l}_{\vec{\omega}} &= \underline{\dot{I}}_t \vec{\omega} + \vec{\omega} \times \left( \underline{I}_t \vec{\omega} + \vec{\kappa}_t + \vec{h}_w \right) + \dot{\vec{\kappa}}_{\ddot{\chi}} \\ \vec{l}_{\vec{\chi}} &= \underline{K}\,\vec{\chi} - \underline{D}\,\dot{\vec{\chi}} + \vec{C} \end{aligned} \tag{27}$$

and $\vec{\omega}^T = [\omega_x\ \omega_y\ \omega_z]$ is the spacecraft angular rate, $\vec{\chi}$ is the non-dimensional generalized displacement coordinates, $\underline{I}_t$ is the spacecraft total inertia tensor, $\underline{\dot{I}}_t$ is the rate of change of the inertia tensor, $\vec{\kappa}_t$ is the angular momentum due to flexibility, $\dot{\vec{\kappa}}_t = \dot{\vec{\kappa}}_{\dot{\chi}} + \dot{\underline{\kappa}}_{\ddot{\chi}}^T \ddot{\vec{\chi}}$ is the rate of angular momentum due to flexibility, $\vec{h}_w$ is the reaction wheels' angular momentum, if any, $\underline{K}$ is the flexible plate appendage stiffness matrix, $\underline{D}$ is the flexible plate appendage's diagonal damping matrix with entries $d$ for the damping parameter, $\vec{C}$ is the flexible appendage's remaining nonlinear terms, $\vec{\tau}$ represents the input control torque, $\underline{M}_{\phi\psi} = ab\underline{I}$, $a$ and $b$ are the width and the length of a rectangular panel (refer to Figure 1 for details), and $\underline{I}$ is the identity matrix.

By defining $\vec{\delta}^T = [\vec{\chi}^T\ \dot{\vec{\chi}}^T]$, (26) can be expressed in a concise state space form (1) as,

$$\begin{bmatrix} \dot{\vec{\omega}} \\ \dot{\vec{\delta}} \end{bmatrix} = \begin{bmatrix} \vec{f}_{\vec{\omega}} \\ \vec{f}_{\vec{\delta}} \end{bmatrix} + \begin{bmatrix} \underline{G}_{\vec{\omega}} \\ \underline{G}_{\vec{\delta}} \end{bmatrix} \vec{\tau} \triangleq \vec{f}(\vec{x}) + \underline{G}(\vec{x})\,\vec{\tau} \tag{28}$$

where the $\vec{f}$ and $\underline{G}$ terms can be obtained directly from (26) and the state and output vectors are given as

$$\vec{x} = [\vec{\omega}^T \quad \vec{\delta}^T]^T = [\omega_x \ \omega_y \ \omega_z \ \chi_1 \ \chi_2 \cdots \chi_{mpq} \ \dot{\chi}_1 \ \dot{\chi}_2 \cdots \dot{\chi}_{mpq}]^T$$
$$\vec{y} = \vec{h}(\vec{x}) = [x_1 \quad x_2 \quad x_3]^T = [\omega_x \quad \omega_y \quad \omega_z]^T \qquad (29)$$

This is a $n$-dimensional system with $p = s = r = 3$ and $n = 3 + 2mpq$, where $m$ is the number of flexible appendages, $p$ and $q$ are the number of discretized modes along $x$ and $y$ appendage axes, respectively.

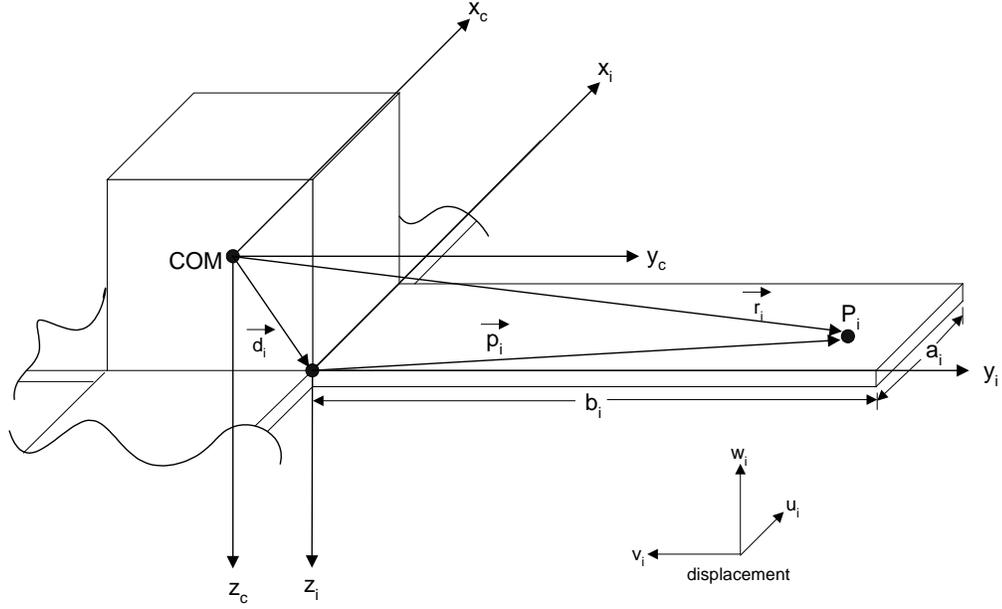

Figure 1.    Spacecraft configuration and appendage reference frame

It can be shown that our system has well-defined relative degrees, $r_i$, for $i = 1, 2, 3$ at $\vec{x}_o$ such that the conditions (30) are satisfied and input-output feedback linearization implementation is possible, namely

(i)   The scalar functions $L_{\vec{g}_j} L_{\vec{f}}^k h_i = 0$ for $\forall i, j, k$ where $1 \leq i, j \leq 3$ and $k \leq r_i - 1$

(ii)   The decoupling matrix $E = \begin{bmatrix} L_{\vec{g}_1} L_{\vec{f}}^{r_1-1} h_1 & L_{\vec{g}_2} L_{\vec{f}}^{r_1-1} h_1 & L_{\vec{g}_3} L_{\vec{f}}^{r_1-1} h_1 \\ L_{\vec{g}_1} L_{\vec{f}}^{r_2-1} h_2 & L_{\vec{g}_2} L_{\vec{f}}^{r_2-1} h_2 & L_{\vec{g}_3} L_{\vec{f}}^{r_2-1} h_2 \\ L_{\vec{g}_1} L_{\vec{f}}^{r_3-1} h_3 & L_{\vec{g}_2} L_{\vec{f}}^{r_3-1} h_3 & L_{\vec{g}_3} L_{\vec{f}}^{r_3-1} h_3 \end{bmatrix}$ is invertible in the neighborhood $\Omega$ of $x_o$ (30) where the Lie derivative $L_{\vec{f}} h$ for the scalar function $h$ along the vector field $\vec{f}$ is define as

$$L_{\vec{f}} h = \nabla h \ \vec{f} = \frac{\partial h}{\partial \vec{x}} \vec{f} = \begin{bmatrix} \frac{\partial h}{\partial x_1} & \cdots & \frac{\partial h}{\partial x_m} \end{bmatrix} \begin{bmatrix} f_1 \\ \vdots \\ f_m \end{bmatrix} = \frac{\partial h}{\partial x_1} f_1 + \cdots + \frac{\partial h}{\partial x_m} f_m$$

We can obtain an input-output feedback linearizing torque for the flexible spacecraft given by [18]

$$\vec{\tau} = \underline{G}_{\vec{\omega}}^{-1}[-\vec{f}_{\vec{\omega}} + \vec{v}] \tag{31}$$

By substituting the control torque (31) into the system dynamics (28), we will obtain (32) which is a simple linear relationship between the output and the new input which constitutes the new external and observable dynamics of the system, and (33) which is an internal unobservable dynamics (i.e. dynamics of $\vec{\chi}, \dot{\vec{\chi}}$ states):

$$\dot{\vec{\omega}} = \vec{v} \tag{32}$$

$$\dot{\vec{\delta}} = \vec{f}_{\vec{\delta}} - \underline{G}_{\vec{\delta}}\,\underline{G}_{\vec{\omega}}^{-1}\vec{f}_{\vec{\omega}} + \underline{G}_{\vec{\delta}}\,\underline{G}_{\vec{\omega}}^{-1}\vec{v} \tag{33}$$

We need to show that both parts of the closed-loop system dynamics are stable together. The proof of stability of (33) will indicate that our closed-loop system is minimum phase and hence there is at least a possibility that the overall system is asymptotically stable. However, before analyzing (33), there is a need to simplify it by finding appropriate coordinate transformations (i.e. diffeomorphisms) to rearrange the system in the "normal form".

The necessary geometric control concepts and the stability proof for flexible spacecraft are presented below:

Using the outputs $\vec{y}$ of the system (28)-(29), their derivatives, the relative degrees $r_i = 2$ for $i = 1, 2, 3$ and the total relative degree of $r = \sum_{i=1}^{3} r_i = 6$, one can devise a non-linear coordinate transformation $T$ to map $\vec{x}$ in (29) to $\vec{\mu}$, as given by

$$\vec{\mu} = T(\vec{x}) = T\begin{pmatrix} q_i \\ \dot{q}_i \\ --- \\ \chi_1 \\ \vdots \\ \dot{\chi}_{npq} \end{pmatrix} = \begin{bmatrix} \Psi_i(\vec{x}) = h_i \\ \Psi_{i+3}(\vec{x}) = L_{\vec{f}}h_i \\ --- \\ \Phi_1(\vec{x}) \\ \vdots \\ \Phi_{m-8}(\vec{x}) \end{bmatrix} = \begin{bmatrix} q_i \\ \dot{q}_i \\ --- \\ \Phi_1(\vec{x}) \\ \vdots \\ \Phi_{m-8}(\vec{x}) \end{bmatrix} \triangleq \begin{bmatrix} \vec{\zeta} \\ --- \\ \vec{\eta} \end{bmatrix} \quad \text{where } i = 1,2,3 \tag{34}$$

with $\vec{\zeta}$ representing the new quaternion-based rigid states (i.e. $\vec{q}$ and $\dot{\vec{q}}$ as opposed to body angular rate $\vec{\omega}$, refer to [19]) and $\vec{\eta}$ representing the remaining states which no longer have the same clear physical meaning as the elastic states $\vec{\chi}$ and $\dot{\vec{\chi}}$. The coordinate transformation $T(\vec{x})$ must be a diffeomorphism which means that $T(\vec{x})$ must to be a continuous, differentiable and bijective (i.e. onto and one-to-one) mapping and the inverse mapping $T^{-1}(\vec{x})$ must also exist. Using our new states, our closed-loop system (32)-(33) can be converted into the following "normal form"

$$\begin{bmatrix} \dot{\vec{\zeta}} \\ --- \\ \dot{\vec{\eta}} \end{bmatrix} \triangleq \begin{bmatrix} \dot{\zeta}_i \\ \dot{\zeta}_{i+3} \\ --- \\ \dot{\vec{\eta}} \end{bmatrix} = \begin{bmatrix} \zeta_{i+3} \\ v_i \\ --- \\ \vec{\alpha}(\vec{\zeta},\vec{\eta}) + \vec{\beta}(\vec{\zeta},\vec{\eta})\vec{v} \end{bmatrix} \quad (35)$$

where the $\dot{\vec{\zeta}}$ and $\dot{\vec{\eta}}$ represent the external and internal dynamics of the closed-loop system. If the distribution $\Gamma$, spanned by smooth vector fields $\vec{g}_j$, $j = 1, 2, 3$ is involutive near $\vec{x}_o$ then a diffeomorphism $\Phi(\vec{x})$ exists such that $\vec{\beta}(\vec{\zeta},\vec{\eta}) = 0$, and therefore the internal dynamics has the simpler form of $\dot{\vec{\eta}} = \vec{\alpha}(\vec{\zeta},\vec{\eta})$ which is now independent of the input $\vec{v}$. The zero dynamics is then defined as [4]

$$\dot{\vec{\eta}} = \vec{\alpha}(0,\vec{\eta}) \quad (36)$$

with output $y_i = q_i = \zeta_i$, $i = 1, 2, 3$ constrained to be equal to zero (also then $\dot{q}_i = \zeta_{i+3} = 0$) for all times by proper choice of initial conditions $\vec{\zeta}(t_o)$ and feedback $\vec{v}(0,\vec{\eta})$.

By definition, an involutive distribution implies that it is closed under the Lie bracket defined as $[\vec{g}_i, \vec{g}_j] = L_{\vec{g}_i}\vec{g}_j - L_{\vec{g}_j}\vec{g}_i$. In our case, it means that the dimension of the distribution $\Gamma = \text{span}\{g_j\}$, $j = 1,2,3$ is always 3 and does not change when a new vector field $\vec{r}$, which is generated by the Lie bracket of any combination of two of the vector fields $\vec{g}_j$, is included in the set of the vector field generating the distribution $\Gamma$. It can be shown that this is in fact the case for our system, in other words, every new Lie bracket generated vector field $\vec{r}$ is always dependent on vector fields $\vec{g}_j$. Given that our distribution $\Gamma$ is involutive, then by the Frobenius Theorem [4], the set of differential equations

$$L_{\vec{g}_j}\Phi_k = \frac{\partial \Phi_k}{\partial \vec{x}}\vec{g}_j = 0, \quad k = 1, 2, \ldots, m-8; \quad j = 1, 2, 3 \quad (37)$$

has a full set of solutions. These solutions are the appropriate diffeomorphisms $\Phi(\vec{x})$, which result in $\vec{\beta}(\vec{\zeta},\vec{\eta}) = 0$. Also note that we further require the following two conditions be met for the appropriate $\Phi(\vec{x})$, namely

(i) $\Phi_k(0) = 0$ $k = 1, 2, \ldots, m-8$, and

(ii) $\nabla \Phi_k$ are linearly independent.

The $n-r$ dimensional zero dynamics of (35) can be directly obtained in a simple and concise closed-form using (25) as

$$\dot{\vec{\eta}} = \underline{F}_{22}^{-1}\vec{l}_{\bar{\chi}} \quad (38)$$

It can be shown that (38) is asymptotically stable, a result that is crucial for the stability analysis of the overall closed-loop system [18].

IV. CONCLUSIONS

The normal form and zero dynamics are powerful tools useful in analysis and control of both linear and nonlinear systems. Different approaches presently found in the literature are all algorithmic and computational in nature and difficult to implement and use for complex high dimensional systems. In contrast in this paper, the author has presented the *first algebraic approach* to this problem where the solutions (24)-(25) are in closed form so that the application of the theorem presented is straight forward for systems with an arbitrary large dimension. An example of the application of the methodology developed underlines the importance of having a simple formulation for obtaining the zero dynamics of highly complex nonlinear and high dimensional systems such as the attitude control of a flexible spacecraft. The zero dynamics (38) would have been extremely difficult to obtain directly from (32)-(33) without the use of the theorem presented in this paper.